\documentclass[a4paper, 10pt]{article}
\usepackage{amssymb}
\usepackage{amsmath}
\usepackage{amscd}
\usepackage{latexsym}
\usepackage{theorem}
\theoremheaderfont{\scshape}
\newtheorem{theorem}{\bfseries Theorem}[section]
\newtheorem{proposition}[theorem]{\bfseries Propositon}

\newtheorem{corollary}[theorem]{\bfseries Corollary}

\theorembodyfont{\rmfamily}
\newtheorem{definition}[theorem]{\bfseries Definition}
\newtheorem{example}[theorem]{Example}

\theoremheaderfont{\rmfamily}
\newtheorem{pf}{Proof}

\numberwithin{equation}{section}

\def\rmod{{R\textrm{-mod}}}
\def\Hom{\mathrm{Hom}}
\def\Ext{\mathrm{Ext}}
\def\Tor{\mathrm{Tor}}
\def\End{\mathrm{End}}
\def\coker{\mathrm{Coker}}

\def\Kdim{\mathrm{dim}}

\def\m{{\mathfrak m}}
\def\n{{\mathfrak n}}
\def\P{\Bbb P}
\def\depth{\mathrm{depth\,}}
\def\gdim{{\mathrm{G}\textrm{-}\mathrm{dim}}}
\def\pd{{\mathrm{pd}}}
\def\CIdim{{\mathrm{CI}\textrm{-}\mathrm{dim}}}

\def\pd{\mathrm{pd}}

\def\rad{{\mathrm{rad}}}

\def\G{{\cal G}}
\def\Z{{\Bbb Z}}
\def\N{{\Bbb N}}
\def\O{{\mathcal{O}}}
\def\qed{$\square$}

\begin{document}
%\begin{article}

%%%%%%%%%%%%%%%%%%%% title %%%%%%%%%%%%%%%%%%%%%
%\begin{opening}

\title{Modules of G-dimension zero over local rings 
 with the cube of maximal ideal being zero
}
\author{Yuji Yoshino\\
Okayama University, 
Faculty of Science, \\
Okayama 700-8530, Japan\\
\texttt{yoshino@math.okayama-u.ac.jp}
}

\maketitle

%%%%%%%%%%%%%%%%%%%%%%%%%%%%%%%%%%%%%%%%%%%%%%%%%
%%%%%%%%%%%%%%%%%%%%%%%%%%%%%%%%%%%%%%%%%%%%%%%%%
%%%%%%%%%%%%%%%%%%%%%% abstract %%%%%%%%%%%%%%%%%
%%%%%%%%%%%%%%%%%%%%%%%%%%%%%%%%%%%%%%%%%%%%%%%%%
%%%%%%%%%%%%%%%%%%%%%%%%%%%%%%%%%%%%%%%%%%%%%%%%%

\begin{abstract}
Let  $(R, \m)$  be a commutative Noetherian local ring with  $\m^3 =(0)$. 
We give a condition for  $R$  to have a non-free module of G-dimension zero.
We shall also construct a family of non-isomorphic indecomposable modules of G-dimension zero  with parameters in an open subset of projective space.
We shall finally show that  the subcategory consisting of modules of G-dimension zero  over  $R$  is not necessarily a contravariantly finite subcategory in the category of finitely generated $R$-modules.  
\end{abstract}

%\end{opening}

%%%%%%%%%%%%%%%%%%%%%%%%%%%%%%%%%%%%%%%%%%%%%%%%%
%%%%%%%%%%%%%%%%%%%%%%%%%%%%%%%%%%%%%%%%%%%%%%%%%
%%%%%%%%%% Introduction %%%%%%%%%%%%%%%%%%%%%%%%%
%%%%%%%%%%%%%%%%%%%%%%%%%%%%%%%%%%%%%%%%%%%%%%%%%
%%%%%%%%%%%%%%%%%%%%%%%%%%%%%%%%%%%%%%%%%%%%%%%%%

\section{Introduction}

The notion of G-dimension of finitely generated modules are introduced by Auslander and Bridger \cite{AB}.
Although various properties concerning G-dimension has been known, 
it still seems to lack references which show us many examples of modules of G-dimension zero. 
One of our aims here is to give such examples of modules.
Actually we construct a family of modules of G-dimension zero with 
continuous parameters for certain cases.

If the ring  $R$  is a Gorenstein local ring, then an $R$-module has G-dimension zero if and only if it is a maximal Cohen-Macaulay module.
Thus it will be natural to expect that many of the properties of the category of maximal Cohen-Macaulay modules over a Gorenstein local ring are satisfied as well for the category of modules of G-dimension zero over a general local ring.
This is a main reason for us to consider modules of G-dimension zero in this paper.
However, as a result of our construction of family of modules G-dimension zero, we conclude unfortunately that, for a certain kind of local ring, the category of modules of G-dimension zero is not a contravariantly finite subcategory in the category of finitely generated modules.

%%%%%%%%%%%%%%
\bigskip
\footnoterule
\medskip
{\footnotesize{
2000 {\it Mathematics Subject Classification}, 
Primary {13C14, 13D05, 16G50}
\par
{\it Key Words}, 
G-dimension, Cohen-Macaulay ring with minimal multiplicity, contravariantly finite subcategory
}}
\newpage
%%%%%%%%%%%%%%

We start recalling the definition and elementary properties of G-dimension in \S 2.
We shall also remark there that if  $(R, \m)$  is an Artinian local ring with  $\m ^2 = (0)$  or more generally if  $R$  is a Cohen-Macaulay local ring with minimal multiplicity, then there is no nonfree $R$-module of G-dimension zero.
 Therefore we have enough reason to consider local rings  $(R, \m)$  with  $\m ^3 =(0)$  as the easiest cases of Artinian local rings which  possess nonfree modules of G-dimension zero.
In \S 3 we get a necessary condition for such $R$ to have nonfree module of G-dimension zero.
Roughly speaking, the ring  is homologically unique, that is,  $R$ should be a Koszul algebra with unique Poincar\'e series, Bass series  and Hilbert series.

In \S 4 we consider mainly Artinian local rings with nontrivial deformation.
For such local rings, the necessary condition that we get in \S 3 is sufficient as well to  have a nonfree module of G-dimension zero.
Surprisingly enough, we can also show that such a local ring  $R$  always  has the form  $S/fS$  where  $S$  is a 1-dimensional Cohen-Macaulay local ring with minimal multiplicity and  $f$  is any nonzero divisor of degree two. 
For the rings of this form, we shall very concretely construct a family of modules of G-dimension zero with continuous parameters in \S 5.
And in \S 6 we prove that the category of modules of G-dimension zero  is not a contavariantly finite subcategory in the category of finitely generated modules over a local ring of this form.

%%%%%%%%%%%%%%%%%%%%%%%%%%%%%%%%%%%%%%%%%%%%%%%%%
%%%%%%%%%%%%%%%%%%%%%%%%%%%%%%%%%%%%%%%%%%%%%%%%%
%%%%%%%%%%%%%% Preliminaries %%%%%%%%%%%%%%%%%%%%
%%%%%%%%%%%%%%%%%%%%%%%%%%%%%%%%%%%%%%%%%%%%%%%%%
%%%%%%%%%%%%%%%%%%%%%%%%%%%%%%%%%%%%%%%%%%%%%%%%%

\section{Preliminaries for G-dimension}

Throughout the present paper, $R$ denotes a commutative Noetherian local ring
 with maximal ideal  $\m$, and by  an $R$-modules we always  mean a finitely generated $R$-module.

In this section, let us recall the definition of G-dimension and its elementary properties.

\begin{definition}\label{def of gdim 0}
{\rm 
Let  $M$ be an $R$-module.
We say that $M$ has G-dimension zero  if it satisfies the following conditions:

\begin{itemize}
\item[(1)] The natural morphism $M \to \Hom_R(\Hom _R(M,R),R)$ is an isomorphism. \vspace{6pt}
\item[(2)] $\Ext ^{i}_R(M,R)=0$ for all $i>0$.\vspace{6pt}
\item[(3)] $\Ext ^{i}_R(\Hom _R(M,R),R)=0$ for all $i>0$. 
\end{itemize}
}\end{definition}

It is easy to see from the definition that if  $M$  is a module of G-dimension zero,   then there is an exact sequence of free $R$-modules;
$$ 
F_{\bullet} : \cdots \to F_n \to F_{n-1} \to \cdots \to F_1 \to F_0 \to F_{-1} \to \cdots 
$$
with the property that the dual complex  $\Hom _R (F_{\bullet}, R)$ is also exact, and  $M \cong  \coker \left( F_1 \to F_0 \right)$.
We call such an exact sequence  $F_{\bullet}$  a complete resolution of  $M$.

\begin{definition}\label{def of gdim}
{\rm 
Let  $n$  be an integer.
We say that an $R$-module  $M$ has G-dimension at most  $n$, denoted by 
$\gdim _R M \leqq n$,   if 
the $n$-th syzygy module  $\Omega ^n _R(M)$  of  $M$  has G-dimension zero.
If there is no such integer  $n$, we denote  $\gdim _R M = \infty$.
}\end{definition}

As it is shown in the paper of Auslander and Bridger \cite{AB}, the G-dimension enjoys  several good properties.
We recall some of them from \cite{AB} for the later use.

\begin{proposition}\label{properties of gdim}
\begin{itemize}
\item[$(1)$]
If  $\gdim _R M < \infty$, then  $\gdim _R M = \depth R - \depth M$.
\vspace{6pt}
\item[$(2)$]
The following are equivalent for a local ring  $R$ : 
\begin{itemize}
\item[$\bullet$] $R$  is a Gorenstein ring.\vspace{6pt}
\item[$\bullet$] $\gdim _R \ M < \infty$  for every $R$-module $M$.\vspace{6pt}
\item[$\bullet$] $\gdim _R \ R/ \m < \infty$.
\end{itemize}
\vspace{6pt}
\item[$(3)$]
Let  $x \in \m$  be a non-zero divisor on both  $R$  and $M$.
Then, we have 
$\gdim _R \ M = \gdim _{R/xR} \ M/xM$.
\vspace{6pt}
\item[$(4)$]
Let  $x \in \m$  be a non-zero divisor on  $R$  satisfying  $xM =0$.
Then, we have 
$\gdim _R \ M = \gdim _{R/xR} \ M + 1$.
\end{itemize}
\end{proposition}

%%%%%%%%%%%%%%%%% In the case of minimal multiplicity %%%%%%%%%%%%%%%%%%%%

Now we assume that  $(R, \m)$ is a local ring that satisfies $\m ^2=(0)$.
In this case, it can be easily seen that there is no nontrivial module of G-dimension zero, unless  $R$  is  Gorenstein.
More precisely, we can prove the following proposition:

\begin{proposition}\label{m^2=0}
Let  $(R, \m)$  be a local ring with  $\m ^2 =(0)$  and suppose that  $R$  is not a Gorenstein ring.
Then every $R$-module of G-dimension zero is a free $R$-module.
\end{proposition}

\vspace{6pt}
\begin{pf}
Let  $M$ be an indecomposable  $R$-module with $\gdim _R M =0$.
Considering the complete resolution of  $M$, we can embed  $M$ into a free $R$-module; $M \subset  F$.
If there is an element  $x$  of $M$ that does not belong to  $\m F$, 
then  $Rx$  is a direct summand of  $F$, hence of  $M$.
Since we assume that  $M$  is indecomposable,  $M \cong  R$  in this case.

Thus we may assume that  $M \subseteqq \m F$.
Then, since  $\m ^2 =(0)$, we have  $\m M =0$, and hence  $M \cong R/\m$, because  $M$  is indecomposable.
As a result, we have  $\gdim  _R R/\m = 0$.
Then it follows from Proposition \ref{properties of gdim}(2) that  $R$  must  
be a Gorenstein ring.
\qed
\end{pf}

Recall that a Cohen-Macaulay local ring  $(R, \m)$  has minimal multiplicity if  there is a system of parameters  $\underline{x} = x_1, x_2, \ldots, x_d$  such that  $\underline{x} \m = \m ^2$.
The above proposition  automatically  implies the same result for modules over 
a Cohen-Macaulay local ring with minimal multiplicity.

\begin{corollary}\label{minimul}
Let  $(R, \m)$  be a Cohen-Macaulay local ring with minimal multiplicity.
Suppose  that  $R$  is not a Gorenstein ring.
Then every $R$-module of G-dimension zero is a free $R$-module.
\end{corollary}

\vspace{6pt}
\begin{pf}
Let  $M$  be an $R$-module with  $\gdim _R M =0$.
Note that  $M$ is a maximal Cohen-Macaulay module over  $R$, for  $\depth M = \depth R$  by Proposition \ref{properties of gdim}(1).
Now take a system of parameters  $\underline{x}$  such as  $\underline{x} \m = \m^2$, and let  $\overline{R}$  (resp. $\overline{M}$) denote $R/\underline{x}R$  (resp. $M/\underline{x}M$).
Note that the square of the maximal ideal of  $\overline{R}$ is zero. 
On the other hand, by a successive use of Proposition \ref{properties of gdim}(3), we  have  $\gdim _{\overline{R}}\overline{M} = 0$.
Thus  $\overline{M}$  is a free $\overline{R}$-module.
Now let $r = \Kdim _k (M \otimes _R k) =  \Kdim _k (\overline{M} \otimes _{\overline{R}} k)$,  where   $k$  denotes the residue class field  $R/\m$.
And we can make a minimal free cover of  $M$ ;  
$$
0 \to N \to R^r \overset{\pi}\to M \to 0,
$$
satisfying  $\pi \otimes _R \overline{R}$  is an isomorphism.
Thus, since  $\Tor_1^R (M, \overline{R}) =0$, we have  $N \otimes _R \overline{R} =0$  and hence  $N = 0$.
Therefore  $M \cong R ^r$  as desired.
\qed
\end{pf}

%%%%%%%%%%%%%%%%%%%%%%%%%%%%%%%%%%%%%%%%%%%%%%%%%
%%%%%%%%%% Section 3 %%%%%%%%%%%%%%%%%%%%%%%%%%%%
%%%%%%%%%%%%%%%%%%%%%%%%%%%%%%%%%%%%%%%%%%%%%%%%%
%%%%%%%%%%%%%%%%%%%%%%%%%%%%%%%%%%%%%%%%%%%%%%%%%
%%%%%%%%% local rings with m^3=0 %%%%%%%%%%%%%%%%
%%%%%%%%%%%%%%%%%%%%%%%%%%%%%%%%%%%%%%%%%%%%%%%%%
%%%%%%%%%%%%%%%%%%%%%%%%%%%%%%%%%%%%%%%%%%%%%%%%%

\section{Local rings with  $\m ^3 =(0)$}

As we have shown in the previous section, Artinian non-Gorenstein local rings with $\m^2 =(0)$  have no nontrivial modules of G-dimension zero.
For the next step, we shall consider Artinian non-Gorenstein local rings with  $\m ^3 = (0)$,  in this section.
If such a ring has a nonfree module of G-dimension zero, then 
the ring structure will be strongly restricted as we prove in the following theorem.

\begin{theorem}\label{exist}
Let  $(R, \m)$  be a non-Gorenstein local ring with  $\m ^3=(0)$, but $\m ^2 \not= (0)$.
For the simplicity we assume that  $R$  contains a field $k$ isomorphic to  $R/\m$.
And let  $r$  be the type of  $R$, i.e. 
$r = \Kdim _k \Hom _R(k, R)$.

Now assume that there is a nonfree $R$-module with  $\gdim _R M=0$.
Then the following conditions hold:

\begin{itemize}
\item[$(1)$]
$R$  has a natural structure of homogeneous graded ring with 
$R = R_0 \oplus R_1 \oplus R_2$, where 
$R_0 =k$ and  $\Kdim _k R_1 = r+1$,  $\Kdim _k R_2 = r$.
In particular, the Hilbert series  of  $R$ is 
$$
H_R (t) = (1+t)(1+rt),
$$
and  $(0:_R \m) = \m^2$.
\vspace{6pt}
\item[$(2)$]
Under the above graded structure,  $R$  is a Koszul algebra
whose Poincar\'e series is 
$$
P_R(t) = \frac{1}{(1-t)(1-rt)}.
$$
\item[$(3)$]
The Bass series of  $R$  is as follows:
$$
B_R (t) = \frac{r-t}{1-rt}
$$
\item[$(4)$]
Every $R$-module  $M$  of G-dimension zero  has a natural structure of graded  $R$-module.
If $M$  has no free summand, then $M$  has only two graded pieces;  
$M = M_0 \oplus M_1$, 
where putting  $b = \Kdim _k M_0$, we see that 
$b$ is equal to the minimal number of generators of  $M$  and 
$\Kdim _k M_1 = rb$, in particular, the length of  the $R$-module $M$  
is  $\ell _R (M) = b(1+r)$.
Furthermore, $M$ has a free resolution of the following type:
$$
\cdots \to R(-n-1) ^b \to R(-n)^b \to \cdots \to  R(-1)^b \to R^b \to M \to 0.
$$
\end{itemize}
\end{theorem}

%%%%%%%%%%%%%%%%%%%%%%% Proof of Theorem 1 %%%%%%%%%%%%%%%%%

\vspace{6pt}
\begin{pf}
Let  $M$  be a nonfree  $R$-module with  $\gdim _R M =0$.
If necessary, taking a summand of  $M$, we assume that 
$M$ is a nonfree indecomposable module.
We prove the theorem step by step.

%%%%%%%%%%%%%%%%%%%%%% step 1 %%%%%%%%%%%%%%
\vspace{12pt}
\noindent
(Step 1) \ We show that  $\m ^2 M =0$.

\vspace{6pt}
In fact, as in the proof of Proposition \ref{m^2=0}, 
$M$  can be embedded into  $\m F$  for some free module  $F$.
Since  $\m^3 =(0)$,  we have  $\m^2 M =0$.
\qed
\vspace{6pt}
 
Therefore, considering a filtration  $(0) \subseteq \m M \subseteq M$,
and noting from the above that  $\m (\m M) =0$, we have the following exact sequence of  $R$-modules:

\begin{equation}\label{filtration of M}
0 \to k^a \to M \to k^b \to 0,
\end{equation}
where  $a = \Kdim _k \m M$ and  $b = \Kdim _k M/\m M$.
Note that the length  $\ell _R(M)$  of the $R$-module  $M$  is  $a+b$.
Note also that  $a \not=0$  and  $b \not= 0$. 
Otherwise,  $M$  would be a direct sum of  $R/\m$  and $R$  would be a Gorenstein ring by Proposition \ref{properties of gdim}(2).

%%%%%%%%%%%%%%%%%%% step 2 %%%%%%%%%%%%%%%%
\vspace{12pt}
\noindent
(Step 2) \ We put 
\begin{equation}\label{def of mu}
\mu _i = \Kdim _k \Ext ^i_R (k, R)  \quad \textrm{for}\quad i \geqq 0.
\end{equation}
And we claim that  
$$
a \geqq b \quad \textrm{and} \quad a\mu _i = b \mu _{i+1} \ \ \textrm{for}  \ \  i \geqq 1.
$$ 

\vspace{6pt}
Before proving this, we note that  $r = \mu _0$.
Since  $\Ext ^i_R(M, R) = 0$  for  $i \geqq 1$, 
it follows from (1) that 
$\Ext ^i_R(k^a , R) \cong \Ext ^{i+1}_R(k^b, R)$  for  $i \geqq 1$, hence 
we have 
$a \mu_i = b \mu _{i+1} \ (i \geqq 1)$.
Thus, 
$$
\mu _1 = 
\left( \frac{b}{a}\right) \mu _2 = 
\left( \frac{b}{a}\right)^2 \mu _3 
= \cdots = 
\left( \frac{b}{a}\right)^i \mu _{i+1}
= \cdots 
$$
Supposed that  $a < b$.
Then we would have  $\mu _1 < \left( \frac{b}{a}\right)^i$  for a large $i$.
Since all the  $\mu _i$  are nonnegative integers, 
we must have  $\mu _i =0$  for all $i \geqq 1$.
This exactly means that  $R$  is a Gorenstein ring, 
which contradicts the assumption.
\qed
\vspace{6pt}

%%%%%%%%%%%%%%%%%% step 3 %%%%%%%%%%%%%%%
\vspace{12pt}
\noindent
(Step 3) \ We denote  $M ^* = \Hom _R (M, R)$.
Then, we have equalities   
$$
\mu _1 = r^2  -1  \quad \textrm{and}\quad   \ell _R(M^*) = \ell_R(M).
$$
\vspace{6pt}	

Applying the functor  $\Hom _R ( \ \ , R)$  to 
the exact sequence  (1), we have an exact sequence 
$0 \to k ^{rb} \to M^* \to k ^{ra} \to k ^{\mu_1b} \to 0$,
 hence 
\begin{equation}\label{filtration of M*}
0 \to k ^{rb} \to M^* \to k ^{ra - \mu_1b} \to 0.
\end{equation}
Note that, since  $M^*$  is also an indecomposable module of G-dimension zero, 
$M^*/\m M^* \cong k ^{ra - \mu_1b}$ and  $\m M^* \cong k ^{rb}$.  
In particular, we have 
\begin{equation}\label{ell tousiki}
\ell (M^*) = r \ \ell (M) - \mu_1 b. 
\end{equation}
Now, since  $\gdim _RM^* =0$, this equality holds for  $M^*$ and hence 
$$
\ell (M^{**})= r \ \ell (M^*) - \mu_1 (r  a - \mu _1 b).
$$
Since  $M^{**} \cong M$, we have 

\begin{eqnarray*}
\ell (M)  &=& r \ \ell (M^*) - \mu_1 (r a - \mu _1 b) \\
&=& r \ \ell (M^*) - \mu_1 (\ell (M^*) - r  b) \\
&=& - \mu_1 \ \ell (M^*) + r ^2 \ \ell (M).
\end{eqnarray*}
As a consequence we get 
\begin{equation}
\mu _1 \ \ell (M^*) = (r ^2 -1) \ \ell (M).
\end{equation}
Since this equation holds for any nonfree indecomposable module  $M$  of G-dimension $0$, one can apply this to  $M^*$ and obtains 
$$
\mu _1 \ \ell (M) = (r ^2 -1) \ \ell (M^*).
$$
Comparing above two equations we finally obtain that $\mu _1 = r ^2 -1$     and  $\ell (M^*) = \ell (M)$, as desired.
\qed

%%%%%%%%%%%%%%%%%% step 4 %%%%%%%%%%%%%%%
\vspace{12pt}
\noindent
(Step 4) \ The following equations hold :
$$
a = r b, \quad \textrm{and} \quad \mu _i = r ^{i-1} \mu_1 = r^{i+1} - r^{i-1} \ \ \textrm{for all} \ \ i \geqq 1.
$$
In particular, the Bass series of  $R$  is 
$$
B_R(t) = \sum _{i \geqq 0} \mu _i t^i 
= r + \sum _{i \geqq 1} (r^{i+1} - r^{i-1})t^i
=  \frac{r-t}{1-rt}
$$

\vspace{6pt}	
 
In fact, we have shown in Step 2 that  
$$
a \mu _i = b \mu _{i+1} \quad \textrm{for} \ \ i \geqq 1,
$$ 
which holds for any nonfree indecomposable module $M$ of G-dimension zero.
Now apply this to  $M^*$, and we have from (3) that 
$$
(r b)\mu _i = (r a -\mu_1 b) \mu_{i+1} 
\quad \textrm{for} \ \ i \geqq 1.
$$
Note that  $\mu _i > 0$  for all $i \geqq 0$, since  $R$  is non-Gorenstein.
See \cite[Lemma (3.5)]{Bass}.
Hence the above two equations implies that 
$(rb)b = a (ra - \mu _1b)$, equivalently 
$
r b^2 + (r^2 -1) ab -ra^2 =0
$
by Step 3.
It follows that 
$(rb -a )(b-ra)=0$.
Since we have shown in Step 2 that  $a \geqq b$, and since $r \geqq 1$, 
this implies that  $a =rb$.
Also from Step 2 we have that 
$\mu _{i+1} = \frac{a}{b} \mu _i = r \mu_{i}$  for  $i \geqq 1$.
Therefore  $\mu _i = r ^{i-1} \mu_1 = r^{i-1}(r^2 -1)$  for  $i \geqq 1$.
\qed

%%%%%%%%%%%%%%%%%% step 5 %%%%%%%%%%%%%%%
\vspace{12pt}
\noindent
(Step 5) \  Since  $M$  is generated by  $b$ elements, we have a short exact 
sequence:
$$
0 \to \Omega M \to R^b \to M \to 0, 
$$
where  $\Omega M $  is the first syzygy module of  $M$.
We claim here that 
$$
\m \Omega M = \m ^2 R^b.
$$
\vspace{6pt}

Since  $R^b \to M$  is a minimal free cover, we have  $\Omega M \subseteqq \m R^b$, therefore  $\m \Omega M \subseteqq \m ^2 R^b$.
On the other hand, since  $\m ^2 M=0$  by Step 1,  we have  $\m ^2 R^b \subseteqq \Omega M$.
In order to prove the above equality, suppose that
$\m \Omega M \not= \m ^2 R^b$.
Then there is an  $x \in \m^2 R^b \backslash \m \Omega M \subseteqq \Omega M$.
Note  that  $x$  is one of the minimal generators of  $\Omega M$  and  $\m x =0$, for  $\m ^3 =(0)$.
Thus the composition of natural maps
$$
k \cong Rx \subseteqq \Omega M \to \Omega M/\m \Omega M \to kx \cong k
$$
is the identity map on  $k$. 
Hence the submodule  $Rx \cong k$  is a direct summand of  $\Omega M$.
It then follows that  $\gdim _R k =0$, hence  $R$  must be a Gorenstein ring.
This is a contradiction, hence we have  
$\m \Omega M = \m ^2 R^b$.
\qed

%%%%%%%%%%%%%%%%%% step 6 %%%%%%%%%%%%%%%
\vspace{12pt}
\noindent
(Step 6) \  Now we prove that $\Omega M$  is minimally generated by $b$ elements and that the following equalities hold:
$$
\Kdim _k \ \m/\m^2 = r+1, \quad \Kdim _k \ \m^2 = r.
$$
\vspace{6pt}	

To prove this, we put  
$e = \Kdim _k \ \m /\m^2$ and  $f = \Kdim _k \ \m ^2$.
First of all, we note from Step 5 that we have the following exact sequence:
\begin{equation}\label{exseq}
0 \to \Omega M / \m \Omega M \to \left( R/ \m^2 \right)^b \to M \to 0
\end{equation}
Therefore, since  $\ell _R (M) = a+b$, it follows that 
$\Omega M / \m \Omega M \cong k ^{be -a}$.
On the other hand, by Step 5, we have  
$\m \Omega M \cong \m ^2 R^b \cong k^{bf}$.
In particular, there is an exact sequence;
$$
0 \to k ^{bf} \to \Omega M \to k^{be-a} \to 0.
$$
Note that  $\Omega M$  is an indecomposable module of G-dimension zero as well as  $M$,  and hence the equations in Step 4 hold for  $\Omega M$.
Thus we have  $bf = r (be-a)$.
Since $a = rb$, we finally get 
$$
f = r(e-r).
$$ 
Note that  $\m ^2 \subseteqq (0:\m)$, since  $\m ^3 =(0)$.
In particular,  $0 < f \leqq r$  from the definition.
Therefore it follows from the above equality that 
$$
f=r  \quad \textrm{and} \quad e-r =1.
$$
Note also that  we have shown that  $\m ^2 = (0: \m)$.
Since  $\Omega M / \m \Omega M \cong k ^{be -a}$,  $\Omega M$ is minimally generated by  $be -a$ elements, but we have 
$be -a = b(r+1)-a = b +(rb -a) = b$, by virtue of Step 4.
\qed

%%%%%%%%%%%%%%%%%% step 7 %%%%%%%%%%%%%%%
\vspace{12pt}
\noindent
(Step 7) \  Now we shall prove that  $R$  has  natural structure of homogeneous graded ring  with  $R = R_0 \oplus R_1 \oplus R_2$, where 
$$
R_0 \cong k, \quad  \Kdim _k R_1 = r+1 \ \ \textrm{and} \ \ 
\Kdim _k R_2 = r.
$$
Furthermore,  $M$  has natural strucure of graded  $R$-module with 
$M = M_0 \oplus M_1$, where 
$\Kdim _k M_0 = b$  and  $\Kdim _k M_1 = rb$.
\vspace{6pt}

Since we have shown in Step 6 that  $\Kdim _k \ \m/\m^2 = r+1$, 
we can find a $k$-linear subspace  $V$  of  $\m$  with  $\Kdim _k V = r+1$ 
in such a way that  $V$  is isomorphic to  $\m /\m^2$  through the natural 
map  $\m \to \m/\m^2$.
Then clearly  $V \cap \m ^2 =(0)$, and hence  
$R  = k \oplus V \oplus \m^2$  as a $k$-vector space.
This gives a graded ring structure of  $R$  with  $R_0 = k$, $R_1 = V$ and  $R_2 = \m^2$.
It is easy to see that  $\m ^2 = V\cdot V$  by the product in  $R$, and thus the graded ring  $R$  is generated in degree $1$, i.e. a homogeneous graded ring.

Similarly to the above, since  $\Kdim _k M/\m M =b$, 
we can find a linear subspace  $U$  of  $M$  with  $\Kdim _k U = b$ 
in such a way that  $U$  is isomorphic to  $M/\m M$  through the natural 
map  $M \to M/\m M$.
Then clearly  $U \cap \m M =(0)$, and hence  
$M = U  \oplus  \m M$ as a $k$-vector space.
It is easy to see that  $V \cdot U = \m M$ and  $V \cdot \m M =0$.
Therefore  $M$ has a graded $R$-module structure  $M_0 \oplus M_1$  with 
$M_0 = U$  and  $M_1 = \m M$.
Note that  $M$  is generated in degree $0$  as a graded  $R$-module.
Note also that  $\Kdim _k \ \m M =a$ and this equals $rb$  by Step 4.
\qed

%%%%%%%%%%%%%%%%%% step 8 %%%%%%%%%%%%%%%
\vspace{12pt}
\noindent
(Step 8) \  We show here that the graded $R$-module $M$  has the following 
type of linear free resolution :
$$
 \cdots \to R(-n-1)^b \to R(-n)^b \to  
 \cdots \to R(-2)^b \to R(-1)^b \to  
 R^b \to M \to 0
$$
In particular, we have an equality 
$$
\Kdim _k \ \Tor _i ^R(M, k)_j  = 
\begin{cases}
b &(i=j), \\
0 &(i \not= j), \\
\end{cases}
$$
where 
$\Tor _i ^R(M, k)_j$  is piece of degree $j$ in the graded module  $\Tor _i ^R(M, k)$.
\vspace{6pt}

We have shown in Step 7 that  $M$  is minimally generated by $b$ elements of  $M_0$.
Therefore we may take a minimal free cover  $R^b \to M$  being a graded homomorphism.
And since 
$0 \to \Omega M \to R^b \to M \to 0$,  $\Omega M$  is also a graded $R$-module.
Furthermore, by the dimension argument, we have 
\begin{eqnarray*}
\Kdim _k (\Omega M) _1 &= & b\ \Kdim _k \ R_1 - \Kdim _k \ M_1 = 
b(r+1)-rb = b,  \\ 
\Kdim _k (\Omega M) _2 &=&  b\ \Kdim _k \ R_2 - \Kdim _k \ M_2 =  rb  \quad  
 \textrm{and}  \\
\Kdim  _k (\Omega M)_i &=& 0 \quad  \textrm{for} \ i \not= 1, 2.
\end{eqnarray*}
Since  $\Omega M$  is also an indecomposable module of G-dimension zero,  it follows from Step 7 that  $\Omega M$  is minimally generated by  $b$ elements in  $(\Omega M)_1$ and that 
there is a minimal free cover  as a graded $R$-module:
$R(-1)^b \to \Omega M$.
Thus one can continue this procedure to get the result stated in Step 8.
\qed

%%%%%%%%%%%%%%%%%% step 9 %%%%%%%%%%%%%%%
\vspace{12pt}
\noindent
(Step 9) \  Now setting  
$\beta _{ij} =  \Kdim _k \ \Tor _i ^R(k, k)_j$
for all integers  $i$  and $j$, we can show that 
\begin{equation}\label{beta}
\beta _{ij} = 
\begin{cases}
1 + r + r^2 + \cdots + r^i  &(i=j \geqq 0), \\
0 &(\textrm{otherwise}). \\
\end{cases}
\end{equation}
In particular,  the graded ring  $R$  is a Koszul algebra whose Poincar\'e series is 
$$
P_R(t) = \frac{1}{(1-t)(1-rt)}.
$$
\vspace{6pt}	

Before proving this, we note from Step 7 that there is an exact sequence 
$0 \to k(-1)^{br} \to M \to k^{b} \to 0$.
Tensoring $k$  over  $R$ to this sequence, we have a long exact sequence of $k$-vector spaces for any integer $i$  and  $j$;  
\begin{equation}\label{long}
\begin{CD}
\cdots \rightarrow \Tor _i^R(k, k)^{br} _{j-1}  \rightarrow  \Tor _i^R(M, k)_j @>{\varphi _{ij}}>> \Tor _i^R(k, k)^b _j  \rightarrow \Tor _{i-1}^R(k, k)^{br} _{j-1}  \rightarrow \cdots  \\
\end{CD}
\end{equation}
Now we shall compute  $\beta _{ij}$.  
First of all, if  $j < 0$,  then it is clear that  $\beta _{ij} =0$  for all $i$.
We prove the equality (7) by induction on $j \geqq 0$.
For the first case we assume  $j =0$.
Putting  $j =0$  in  (8), we see that  $\varphi _{i0}$  is an isomorphism, since  $\beta _{i,-1} = 0$  for all $i$.
Therefore we have 
$ b \beta _{i0} = \Kdim _k \ \Tor _i^R(M, k)_0$, hence it follows from Step 8 that  $\beta _{00} = 1$  and  $\beta _{i0} = 0$  if $i >0$.
  
Now assume  $j >0$. 
By the induction hypothesis, we have that 
$\Tor _i^R(k , k)_{j-1}=0$  unless $i = j-1$.
Hence it follows from (8) that 
\begin{equation*}\label{long2}
\begin{CD}
0 @>>> \Tor _j^R(M, k)_j @>{\varphi _{jj}}>> \Tor _j^R(k, k)^b _j @>>> \Tor _{j-1}^R(k, k)^{br} _{j-1}  \\
@>>>  \Tor _{j-1}^R(M, k)_j @>{\varphi _{j-1 j}}>> \Tor _{j-1}^R(k, k)^b _j @>>> 0,  \\
\end{CD}
\end{equation*}
and that  $\Tor _i^R(M, k)_j \cong \Tor _i^R(k, k)^b _j$  if $i \not= j, \ j-1$.
Using the result in Step 8, we get 
$\Tor _i^R(k, k)^b _j = 0$  if $i \not= j$  and 
$b \beta _{jj} = b + br \beta _{j-1 j-1}$.
This shows that 
$\beta _{ij} =0$  for $i \not= j$  and  $\beta _{jj} = 1+r \beta_{j-1 j-1}$, which proves Step 9.
\qed

\vspace{6pt}
From all these steps we have proved any of the statements in Theorem (\ref{exist}).
(Q.E.D.)
\end{pf}

%%%%%%%%%%%%%%%%%%%%%%%%%%%%%%%%%%%%%%%%%%%%%%%%%
%%%%%%%%%% Section 4 %%%%%%%%%%%%%%%%%%%%%%%%%%%%
%%%%%%%%%%%%%%%%%%%%%%%%%%%%%%%%%%%%%%%%%%%%%%%%%
%%%%%%%%%%%%%%%%%%%%%%%%%%%%%%%%%%%%%%%%%%%%%%%%%
%%%% local rings with nontrivial deformation  %%%
%%%%%%%%%%%%%%%%%%%%%%%%%%%%%%%%%%%%%%%%%%%%%%%%%
%%%%%%%%%%%%%%%%%%%%%%%%%%%%%%%%%%%%%%%%%%%%%%%%%

\section{Local rings with nontrivial deformation}

As we have shown in the previous section, a local ring  $(R, \m)$  with  $\m ^3=(0)$  has a nonfree module of G-dimension zero only if  $R$  is a homogeneous graded ring.
Therefore, in this section,  we consider only a homogeneous graded Artinian ring, that is,  the ring  $R$  is a residue ring of a polynomial ring by a homogeneous ideal; 
$$
R = k 	[ X_0, X_1, \ldots , X_r]/I, 
$$
where  all variables $X_i$  are of degree one and  $I$  is a homogeneous ideal, and  $\m$  is the maximal homogeneous ideal generated by $X_0, X_1.\ldots , X_r$.

In this case, we make the following definition.

\begin{definition}{\rm
we say that a graded Artinian ring  $R$  has a {\bf nontrivial deformation} if 
there is a homogenous Cohen-Macaulay graded ring  $S$  of positive dimension  $d$   and a homogeneous regular sequence  $f_1, f_2, \ldots , f_d$  of  $S$  with  degree $\geqq 2$  such that $R \cong S/(f_1, f_2, \ldots , f_d)S$.
}\end{definition}

In Theorem (\ref{exist}) we have shown several necessary conditions for  $R$  to have a nonfree module of G-dimension zero.
If  $R$  has a nontrivial deformation,  then we can prove they are    sufficient as well.

\begin{theorem}\label{with deformation}
Let  $R$  be a homogeneous graded ring over a field $k$ with the graded maximal  ideal  $\m$.
Assume that  $R$  is not a Gorenstein ring.
And suppose that  $R$  has a nontrivial deformation.
Then the following conditions are equivalent:
\begin{itemize}
\item[$(1)$]
$R$  has a nonfree module of G-dimension zero, and  $\m ^3 =(0)$.
\vspace{6pt}
\item[$(2)$]
The Hilbert series of  $R$  is 
$$
H_R (t) = (1+t)(1+rt),
$$
for some integer  $r \geqq 2$.
\vspace{6pt}
\item[$(3)$]
There is a one-dimensional Cohen-Macaulay homogeneous graded ring  $S$  with minimal multiplicity and a homogeneous non-zero divisor  $f$  of degree $2$  of  $S$  such that 
$R \cong S/fS$.
\end{itemize}
Furthermore, under the condition $(3)$, we have an equality 
\begin{equation}\label{formula} 
\gdim _R \ M  = \pd _S \ M -1 \ \ (\leqq \infty)
\end{equation}
for any $R$-module $M$.
$($ Therefore  $\gdim _R \ M = \CIdim _R \ M$. $)$
\end{theorem}
\vspace{6pt}
\begin{pf}
$(1) \Rightarrow (2)$ : 
We have proved in Theorem (\ref{exist}) the condition on the Hilbert series of $R$. 
We just note that  the socle dimension  $r = \Kdim _k \Hom _R(k, R)$  is greater than $1$, since  $R$  is not Gorenstein.

$(2) \Rightarrow (3)$ : 
Let  $S$  be a Cohen-Macaulay homogeneous graded algebra over  $k$  of positive dimension  $d$  and  
 $R \cong S/(f_1, f_2 , \ldots , f_d)S$  with a regular sequence  $f_1, f_2, \ldots , f_d$  in  $S$  of degree $\geqq 2$.
Note, in general, that the Hilbert series of  $S$  is of the form 
$$
H_S (t) = \frac{p(t)}{(1-t)^d}, 
$$
for some polynomial  $p(t)$  with nonnegative integral coefficients.
See  \cite[Corollary 4.1.10]{BH}.
And also note that  $S$  has minimal multiplicity if and only if  $p(t) = 1 + nt$  for some integer  $n$.  
Putting  $c_i = deg (f_i) \geqq 2$, we know that there is an equality  
$$
H_R (t) = \prod _{i=1}^d (1-t^{c_i}) \ H_S(t).
$$
Hence, under the condition of $(2)$,  we have 
$$
(1+t)(1+rt) = \prod_{i=1}^d (1+t+t^2 + \cdots +t^{c_i-1}) \ p(t).
$$  
Since this equality holds as an elements in  $\Z [t]$, we must have 
$d =1$, $c_1=2$  and  $p(t) = 1+rt$. 
This exactly means that  $S$  is of dimension one, having minimal multiplicity 
and  $f_1$  is of degree $2$.

$(3) \Rightarrow (1)$ : 
Let  $R \cong S/fS$, where $S$ is a one-dimensional Cohen-Macaulay homogeneous graded ring  with minimal multiplicity and $f$  is a homogeneous non-zero divisor of degree $2$  of  $S$.

First we shall prove the equality (9).
For this, let  $M$ be an $R$-module.
Then it follows from Proposition \ref{properties of gdim}(4) that 
$\gdim _R \ M = \gdim _S \ M -1$.
In particular, 
$\gdim _R \ M < \infty$  if and only if  $\gdim _S \ M < \infty$.
If this is the case, then a certain syzygy module of  $M$  as an  $S$-module has G-dimension zero, however as we showed in Corollary \ref{minimul}, such a syzygy module has to be a free $S$-module.
This shows that   $\gdim _R \ M < \infty$  if and only if the $S$-module $M$ has finite projective dimension, i.e. $\pd _S \ M < \infty$.
On the other hand, we know from  Proposition \ref{properties of gdim}(1) that 
 $\gdim _R \ M = \depth R - \depth M = 0$  and 
$\pd _S \ M = \depth S - \depth M =1$  if they are finite.
Thus we have shown the equality (9) including the case that the both sides are infinite.
(Note from \cite{AGP} that it is known the inequality  
$\gdim _R \ M \leqq \CIdim _R \ M \leqq \pd _S \ M -1$  for any $R$-module  $M$, hence they are all equal in this case.)

Next we note that  $\m ^3 =(0)$.
In fact, since $S$  has minimal multiplicity of dimension one, the Hilbert series  $H_S (t)$  of  $S$  is of the form  
$$
H_S (t) = \frac{1+rt}{1-t}.
$$ 
Thus, 
$$
H_R (t) = H_S(t) (1-t^2)= (1+t)(1+rt),
$$  
since  $f$  is a non-zero divisor of degree two.
This shows that  $R$  has no nontrivial homogeneous components of degree greater than two, and hence we have  $\m ^3 =(0)$.

Now we want to show that there exists a nonfree $R$-module of G-dimension zero.
If  $k$  is an algebraically closed field, then we shall construct a continuous family of such modules in the next section.
And the same proof shows the existence of such a module when  $k$  is an infinite field.
Here we just prove the existence of such a module by another idea using matrix factorizations, which can be applied to any case including the case that  $k$  is a finite field.

Since  $f \in S$  is an element of degree $2$, we can write 
$$
f = \sum _{i=1} ^n \ x_i y_i, 
$$
where  $x_i, y_i$  are elements of $S$  of degree $1$.
We set graded free S-modules  $F$  and  $\bigwedge F$  as follows:
$$
F = \bigoplus _{i=1}^n S e_i , \quad  \bigwedge F = \bigoplus _{i=1}^n \bigwedge ^i F
$$
And we define homogeneous $S$-linear maps $\phi$  and  $\psi : \bigwedge F \to \bigwedge F$  by
$$\begin{cases}
&\phi (e_{i_1} \wedge \cdots \wedge e_{i_s}) = 
\sum _{j=1}^s x_{i_j} (e_{i_1} \wedge \cdots \wedge e_{i_{j-1}} \wedge e_{i_{j+1}} \wedge \cdots \wedge e_{i_s}),  \vspace{12pt} \\
&\psi (\omega ) = (\sum _{i=1} ^n y_i e_i) \wedge \omega. \\
\end{cases}
$$
Then it is easy to verify the equalities 
$\phi ^2 = \psi ^2 = 0$ and  $\phi \cdot \psi + \psi \cdot \phi = f \cdot 1_{\bigwedge F}$, hence  
$(\phi +\psi)(\phi + \psi ) = f\cdot 1 _{\bigwedge F}$.
This means that the pair of linear maps  
$(\phi+\psi, \phi+\psi)$  on  $\bigwedge F$  gives  a matrix factorization of  $f$.
See \cite[p.66]{Y}.
Hence, if we define an $S$-module  $M$  by the exact sequence :
$$
\begin{CD}
0 @>>> \bigwedge F @>{\phi + \psi}>> \bigwedge F @>>> M @>>> 0, 
\end{CD}
$$
it is easy to see that  $f$ annihilates  $M$,  hence $M$  is actually an $R$-module and  $\pd _S \ M = 1$.
Hence  $\gdim _R \ M =0$  by (9).
On the other hand, the Hilbert series of  $M$ can be computed as
$$
H_M(t) = H_{\bigwedge F} (t) (1-t)= 
2^n H_S (t) (1-t) = 2^n (1+rt),
$$
and thus  $M$ has only two graded pieces, hence  $M$ never be a free $R$-module.\qed
\end{pf}

\begin{example}{\rm
For any integer  $r \geqq 2$, we can construct examples of rings $R$  satisfying the conditions in  Theorem \ref{with deformation}.
In fact, we set 
$$
S = k[X_0, X_1, \ldots , X_r]/I_2 
\left(
\begin{array}{llllll}
X_0 &X_1 &X_2 &\cdots &X_{r-1} &X_r \\
X_1 &X_2 &X_3 &\cdots &X_{r} &X_0 \\
\end{array}
\right),
$$
where  $I_2 (A)$  denotes the ideal generated by all the 2-minors of the matrix  $A$. Then it is easy to see that  $S$  is a Cohen-Macaulay ring of dimension one with the Hilbert series 
$$
H_S (t) = \frac{1+rt}{1-t},
$$
hence  $S$  has minimal multiplicity.
To get an example, just set  $R= S/fS$  for any homogeneous element $f$  of  $S$  of degree $2$.
}\end{example}

\begin{example}{\rm
In her paper \cite{V}, Veliche considers a ring 
$$
R = k[x, y, z, w]/(x^2, xy-zw, xy-w^2, xz-yw, xw-y^2, xw-yz, xw-z^2),
$$
and an $R$-module defined by the following exact sequence
$$
\begin{CD}
\cdots 
\longrightarrow R^2 @>{
\left(
\begin{array}{ll}
z & x \\
w & y \\
\end{array}
\right)
}>>
R^2 @>{
\left(
\begin{array}{rr}
y & -x \\
-w & z \\
\end{array}
\right)
}>> R^2 @>{
\left(
\begin{array}{ll}
z & x \\
w & y \\
\end{array}
\right)
}>> R^2 \longrightarrow  M \longrightarrow  0. 
\end{CD}
$$
She proved that  $\gdim _R \ M =0$,  but  $\CIdim _R\ M = \infty$.
Note that  $\m ^3 =(0)$  for the graded maximal ideal  $\m$  of  $R$, hence 
$R$  is a Koszul algebra by Theorem \ref{exist}.
However  we also conclude from Theorem \ref{with deformation} that the ring  $R$  has no nontrivial deformations.
}\end{example}

%%%%%%%%%%%%%%%%%%%%%%%%%%%%%%%%%%%%%%%%%%%%%%%%%
%%%%%%%%%% Section 5 %%%%%%%%%%%%%%%%%%%%%%%%%%%%
%%%%%%%%%%%%%%%%%%%%%%%%%%%%%%%%%%%%%%%%%%%%%%%%%
%%%%%%%%%%%%%%%%%%%%%%%%%%%%%%%%%%%%%%%%%%%%%%%%%
%%%% family of modules of G-dimension 0  %%%%%%%%
%%%%%%%%%%%%%%%%%%%%%%%%%%%%%%%%%%%%%%%%%%%%%%%%%
%%%%%%%%%%%%%%%%%%%%%%%%%%%%%%%%%%%%%%%%%%%%%%%%%

\section{Construction of continuous family of modules of G-dimension zero}
In this section, we consider the case that the conditions in  Theorem \ref{with deformation} hold, i.e. we always denotes that 
$S$  is a homogeneous Cohen-Macaulay graded ring over a field  $k$  with minimal multiplicity of dimension one,  $f$  is an homogeneous non-zero divisor of  $S$  of degree two and that  $R \cong  S/fS$.
Recall that the Hilbert series of  $R$  is given as 
$$
H_R (t) = (1+t)(1+rt),
$$
where  $r$  is the socle dimension of  $R$. 

In such a case, we have shown in the previous section that there exists a nonfree $R$-module of G-dimension zero. 

The purpose of this section is to show that under a further condition that   $k$  is an algebraically closed field, we can actually construct a continuous family of isomorphism classes of indecomposable modules of  G-dimension zero. 
The main result is the following.

%%%%%%%%%%%%%%%%Theorem [family]%%%%%%%%%%%%%%%%%%%%%%%%%%%%%%%%%%

\begin{theorem}\label{family}
Let  $R$  be as above with  $k$   being an algebraically closed field  and let  $V$  be the $k$-vector space that consists of elements of degree one of  $R$.
Note that  $\Kdim _k V = r+1$.
Now denote by  $\P (V)$  the projective space over  $V$, i.e. the set of $k$-linear subspaces of dimension one in  $V$. 
Under these circumstances, there are a non-empty open subset $\O$  of  $\P (V)$  and  $R$-modules  $M(p, n)$  for each  $(p, n) \in \O \times \N$  satisfying the following conditions:
\begin{itemize}
\item[$(1)$]
Each  $M(p, n)$  is an indecomposable  $R$-module that is minimally generated by $n$ elements.
\vspace{6pt}
\item[$(2)$]
$\gdim _R \ M(p, n) = 0$  for any  $(p, n) \in \O \times \N$. 
\vspace{6pt}
\item[$(3)$]
If $(p, n) \not=  (p' , n') \in \O \times \N$, then  $M(p, n)$  and  $M(p', n')$  are non-isomorphic.
\end{itemize}
\end{theorem}

\vspace{6pt}

%%%%%%%%%%%%%%%  Proof of Them [family]  %%%%%%%%%%%%%%%%
\begin{pf}
%%%%%%%%%%%%%%%%%% step 1 %%%%%%%%%%%%%%%
(Step 1) \  
First we note that there is a homogeneous element of degree one  $x \in V$  with   $x \n = \n ^2$  where  $\n$  is the graded maximal ideal of  $S$, since  $S$  is a one-dimensional  Cohen-Macaulay ring with minimal multiplicity and since  $k$  is an infinite field.
Now define the subset  $\O_1$  of  $\P (V)$  as 
$$
\O _1 = \{ [x] \in \P (V) \ | \ x \in V, \ x \n = \n ^2 \},
$$
where  $[x]$  denotes the $k$-linear subspace of  $V$  generated by  $x$.
We claim that  $\O_1$  is a non-empty open subset of  $\P (V)$.

\vspace{12pt}
In fact, the multiplication of  $S$  induces a $k$-linear map
$
\phi : V \to \Hom _k (V, S_2).
$   
Fixing a $k$-basis  $\{ e_0, \ldots , e_r \}$ of  $V$  and writing  
$x = \sum _{i=0} ^r a_i e_i \ (a_i \in k)$, we see that 
$\phi (x) : V \to S_2$  is a surjective map if and only if  
the rank of the matrix  $\sum _{i=0}^r a_i \phi (e_i)$  
is not smaller than  $\Kdim _k S_2$.
Thus it is easy to see that the set consisting of elements  $x \in V$  
with  $\phi (x)$  being surjective  is an open subset of  $V$.
Note that an element  $[x] \in \P (V)$  belongs to  $\O_1$  if and only if 
 $\phi (x)$  is surjective,  and that  there is at least one such  $x$  as remarked above.
Therefore  $\O _1$  is a non-empty open subset of  $\P (V)$.
Compare with the proof of \cite[Theorem 14.14]{M}.
\qed

%%%%%%%%%%%%%%%%%% step 2 %%%%%%%%%%%%%%%
\vspace{12pt}
\noindent
(Step 2) \ 
For  $[x] \in \O_1$,  define an  $S$-module  $M([x], 1)$  as 
$$
M([x], 1) = S/xS.
$$ 
Then,  $M([x], 1)$  is actually an $R$-module and  $\gdim _R \ M([x], 1)= 0$.

\vspace{12pt}
Note that any  $x \in \O _1$  is a non-zero divisor on  $S$, since 
it is a parameter of the Cohen-Macaulay ring  $S$.
Thus we have the following exact sequence of  $S$-modules:
$$
\begin{CD}
0 @>>> S(-1) @>x>>  S @>>> M([x], 1) @>>> 0.
\end{CD}
$$
In particular, the projective dimension of  $M([x], 1)$  as an $S$-module is one.
Since  $x \n = \n^2$, and  since  $f$  is of degree two, we have 
that  $f \in \n ^2 \subseteqq xS$, hence that  $f M([x], 1) =0$.
This implies that  $M([x], 1)$  is actually an $R$-module.
Thus it follows from (9)  in  Theorem \ref{with deformation} that 
$\gdim _R \ M([x], 1) =0$.
\qed

\vspace{12pt}
Note that  $M([x], 1)$  is indecomposable, because it is generated by a single element.
Also note that  $M ([x], 1) = R/xR$.

Now we fix a nonzero element  $z \in V$  and we set 
$$
\O = \O _1 \backslash \{ [z] \},
$$
which is also a non-empty open subset of  $\P (V)$.
For an element  $([x], n) \in \O \times \N$, we define an $R$-module  $M([x], n)$  by the following exact sequence:
$$
\begin{CD}
R(-1)^n @>{\Phi _x}>> R^n @>>> M([x], n) @>>> 0, 
\end{CD}
$$
where  $\Phi _x$  is a $n \times n$ matrix 
\begin{equation}\label{def of M}
\font\b=cmr10 scaled \magstep4
\begin{pmatrix}
x & z  & &  &  & \\
 & x  &z &  & \smash{\hbox{\b 0}}& \\
 &   &x & \ddots  & & \\
 &   & &\ddots  & z & \\
 &\smash{\hbox{\b 0}}  & &  &x  & z \\
 &   & &    &  & x \\
\end{pmatrix}.
\end{equation}

%%%%%%%%%%%%%%%%%% step 3 %%%%%%%%%%%%%%%
\vspace{12pt}
\noindent
(Step 3) \ We shall prove that  $\gdim _R \ M([x], n) = 0$  for all  $([x], n) \in \O \times \N$.

\vspace{12pt}
From the definition, we have a filtration of  $M([x], n)$  by submodules;
$$
M_0 = 0 \subset M_1 \subset \cdots \subset M_{n-1} \subset M_n = M([x], n),
$$
such that  $M_i/M_{i-1} \cong R/xR$  for  $1 \leqq i \leqq n$.
Since  $\gdim _R \ R/xR = 0$, we have  $\gdim _R \ M([x], n)=0$ as desired.
\qed

%%%%%%%%%%%%%%%%%% step 4 %%%%%%%%%%%%%%%
\vspace{12pt}
\noindent
(Step 4) \ Next we show that the modules  $M([x], n)$  are indecomposable.

\vspace{12pt}
Let  $\Lambda = \End _R (M([x], n))$  be the endomorphism ring  of $M([x], n)$. Note that  $\Lambda$  is a nonnegatively graded (noncommutative) finitely  dimensional $k$-algebra.
It is enough to prove that  $\Lambda$  is a local ring.
For this, let  $\varphi$  be an element of  $\Lambda$ of degree $0$.
Then  $\varphi$   induces a commutative diagram 
$$
\begin{CD}
R(-1)^n @>{\Phi _x}>> R^n @>>> M([x], n) @>>> 0 \\
@VQVV  @VPVV @V{\varphi}VV \\
R(-1)^n @>{\Phi _x}>> R^n @>>> M([x], n) @>>> 0, 
\end{CD}
$$
where  $P$  and  $Q$  are $n \times n$ square matrices whose entries are in  $k$,  for  $\varphi$  is of degree $0$.
By the commutativity of the diagram, we have
$\Phi _x Q = P \Phi _x$. 
Putting 
\begin{equation}\label{def of EJ}
E = \begin{pmatrix}
1 & & & & \\
  & 1 & & & \\
  & &\ddots  & & \\
 & & & &1 \\
\end{pmatrix} \quad \textrm{and} \quad
J = \begin{pmatrix}
0 &1 & & & \\
  & 0 &1 & & \\
  & &\ddots  &\ddots & \\
 & & &\ddots &1 \\
 & & & &0 \\
\end{pmatrix}, 
\end{equation}
we have 
$xQ + z JQ = x P + z PJ$.
Since  $x$  and  $z$  are linearly independent as elements
 of the $k$-vector space  $V$, it follows that 
$Q = P$  and  $JQ = PJ$, hence  $JP = PJ$. 
Then it is easy to see that 
\begin{equation}\label{P}
P = a_0 E + a_1 J +a_2 J^2 + \cdots + a_{n-1} J^{n-1},
\end{equation}
  for some  $a_i \in k$.
Thus all the elements in  $\Lambda$  of degree $0$ are induces by matrices of the form (12).
Note that  if  $i>0$, then  $J^i$  induces a nilpotent element of  $\Lambda$, hence an element in the radical  $\rad(\Lambda)$.
Since any elements of $\Lambda$ of positive degree are nilpotent as well, they are also belonging to  $\rad (\Lambda)$.
Thus we conclude that  $\Lambda / \rad (\Lambda) \cong k$  whose nontrivial element is represented by  $a_0 E$  in  (12).
As a  result,  $\Lambda$  is a local ring.
\qed

%%%%%%%%%%%%%%%%%% step 5 %%%%%%%%%%%%%%%
\vspace{12pt}
\noindent
(Step 5) \ Finally we show that if  $([x], n) \not= ([x'], n') \in \O \times \N$, then  $M([x], n) \not\cong M([x'], n')$. 
This will complete the proof of the theorem.

\vspace{12pt}
Since  $n = \Kdim _k (M([x], n) \otimes _R R/\m)$, the above claim will be obvious if  $n \not= n'$.
Thus we may assume that  $n = n'$.

Recall, in general, if an $R$-module $M$  has a presentation 
$$
\begin{CD}
R^n @>A>> R^m @>>> M @>>> 0
\end{CD}
$$
with a matrix  $A$, then the ideal generated by all the entries of  $A$  is an invariant of the isomorphism class of the module $M$.
That is actually an invariant called the Fitting invariant.

Now suppose  $M([x], n) \cong M([x'], n)$  for  $[x], [x'] \in \O$.
Then by this remark we must have  $(x, z)R = (x' ,z)R$.
Looking at the degree one part of this, we have 
$[x, z] = [x' , z]$  as a two-dimensional subspace of  $V$.
Thus we can write 
$$
x' = \alpha x + \beta z \quad (\alpha \not= 0 \in k, \ \ \beta \in k).
$$ 
As in the same way as the proof of  Step 4, we have 
 $P$  and $Q \in GL (n, k)$  such that 
$\Phi _x Q = P \Phi _{x'}$, since  $M([x], n) \cong M([x'], n)$.
Consequently, it follows that 
$$x Q + z JQ =  x' P + z PJ
= \alpha x P + \beta z P + z PJ,
$$
where  $J$  is the matrix defined in (11).
Since  $x$  and  $z$   are linearly independent, we have 
$Q = \alpha P$  and  $JQ = \beta P + PJ$.
Hence,  $\alpha P^{-1}JP =  \beta E + J$. 
Note that the left hand side is nilpotent and that the right hand side is of the Jordan standard form with eigen value $\beta$.
Therefore we get  $\beta =0$, hence  $x' = \alpha x$.
This implies  $[x] = [x']$  as an element of  $\O$.
\qed
\end{pf}

%%%%%%%%%%%%%%%%%%%%%%%%%%%%%%%%%%%%%%%%%%%%%%%%%
%%%%%%%%%% Section 6 %%%%%%%%%%%%%%%%%%%%%%%%%%%%
%%%%%%%%%%%%%%%%%%%%%%%%%%%%%%%%%%%%%%%%%%%%%%%%%
%%%%%%%%%%%%%%%%%%%%%%%%%%%%%%%%%%%%%%%%%%%%%%%%%
%%%% NOT contravarintly finite           %%%%%%%%
%%%%%%%%%%%%%%%%%%%%%%%%%%%%%%%%%%%%%%%%%%%%%%%%%
%%%%%%%%%%%%%%%%%%%%%%%%%%%%%%%%%%%%%%%%%%%%%%%%%

\section{Modules of G-dimension zero are not contravariantly finite}
For a  general local ring   $(R, \m)$, we  denote  by 
 $\rmod$  the category of finitely generated $R$-modules and 
we define a full subcategory  $\G (R)$  of $\rmod$  as follows:
$$
\G (R) = \{ M \in \rmod \ | \ \gdim _R \ M = 0 \}.
$$
According to Auslander, we call  $\G (R)$  a contravariantly finite subcategory of  $\rmod$  if it satisfies the following condition:

\vspace{6pt}
\noindent 
(*) \ 
For any  $M \in \rmod$, there is an exact sequence 
\begin{equation}\label{appr}
\begin{CD}
0 @>>> Y @>>> X @>{\pi}>> M @>>> 0,
\end{CD}
\end{equation}
with  $X  \in \G (R)$  such  that any morphism from  any object of $\G (R)$ to $M$ factors through the morphism   $\pi$.
\vspace{6pt}

Note that Wakamatsu's lemma  \cite{W} or \cite[Lemma 2.1.1]{X} says that the condition (*) is satisfied for  a sequence of the form (13) if and only if  $\Ext ^1_R (X', Y) =0$  for any  $X' \in \G (R)$.

The morphism  $\pi : X \to M$  with the property  (*)  is often called  a (right) $\G (R)$-approximation of  $M$. 
 When  $R$  is a Gorenstein ring, then  $\G (R)$  coincides with the full subcategory consisting of all maximal Cohen-Macaulay modules, and it is known by  \cite{AB2} that  $\G (R)$  is contravariantly finite.
Therefore  any module has a $\G (R)$-approximation if  $R$  is Gorenstein.

In this section we claim that  $\G (R)$  may not be a contravariantly finite subcategory of  $\rmod$.

\begin{theorem}
Let  $S$  be a one-dimensional Cohen-Macaulay homogeneous graded ring over a field  $k$, which is non Gorenstein and has  minimal multiplicity.
Taking a minimal reduction  $x$  of the maximal ideal $\n$  of  $S$, i.e.  $x \n = \n^2$,  we set  $R = S/x^2 S$.
Then the $R$-module  $k$  has no $\G (R)$-approximation.
\end{theorem}

\vspace{6pt}
\begin{pf}
To prove the theorem, supposed that there were an exact sequence:
\begin{equation}\label{apprk}
\begin{CD}
0 @>>> Y @>{i}>> X @>{\pi}>> k @>>> 0,
\end{CD}
\end{equation}
with the property  (*).
We may take such a sequence minimal, i.e. there are no common direct summands   of  $Y$  and  $X$  through  $i$. 
Now we decompose  $X$  into indecomposable modules ; 
$X =  \bigoplus _{i=1} ^n X^{(i)}$. 
Note that each  $\pi |_{X^{(i)}} : X^{(i)} \to k$  is nontrivial, because of the minimality of  (14).
Recall from Theorem \ref{exist} that all the modules  $X_i$'s are graded modules, hence  it follows that  $\pi$  is a graded homomorphism of degree $0$.
Now  assume that  $X^{(i)} \cong R$  for  $1 \leqq i \leqq u$  and that 
$X^{(j)}$  are not free for  $u+1 \leqq j \leqq n$.
Then from the results of  Theorem \ref{exist} we have the graded structure  $X_0 \oplus X_1 \oplus X_2$  on  $X$  with 
$$
\Kdim _k \ X_0 = u + \sum _{j=u+1} ^n s_j, \ \Kdim _k \ X_1 = u(r+1) + \sum _{j=u+1} ^n  r s_j, \ \Kdim _k \ X_2 = ru,  
$$ 
where  $s_j$  is the minimal number of generators of  $X^{(j)}$.
Now since  $\pi$  is graded,  $Y$  has also a graded structure;   
$Y = Y_0 \oplus Y_1 \oplus Y_2$  and it follows from  (14) that  
\begin{equation}\label{dim of Y}
\Kdim _k \ Y_0 = u + \sum _{j=u+1} ^n s_j -1, \ \Kdim _k \ Y_1 = u(r+1) + \sum _{j=u+1} ^n  r s_j, \ \Kdim _k \ Y_2 = ru.  
\end{equation}
As in the proof of Theorem \ref{family},  $R/xR$  is an $R$-module of G-dimension zero.	
In fact, it has a  complete resolution of the form:
$$\begin{CD}
\cdots @>x>>  R  @>x>>  R  @>x>>  R @>>> \cdots.
\end{CD}$$
Hence if follows from Wakamatsu's lemma that 
$\Ext ^1 _R (R/xR, Y) =0$, and that the sequence 
\begin{equation*}
\begin{CD}
\cdots @>x>>  Y  @>x>>  Y  @>x>>  Y @>>> \cdots
\end{CD}
\end{equation*}
is an exact sequence of graded  $R$-modules. 
Taking a graded piece of this, 
we have an exact sequence  of $k$-vector spaces:
\begin{equation*}
\begin{CD}
0  @>>>  Y_0  @>x>>  Y_1  @>x>>  Y_2 @>>> 0 
\end{CD}
\end{equation*}
Thus it follows from  (15) that  
$$
\left( u + \sum _{j=u+1} ^n s_j -1 \right) + ru = 
 u(r+1) + \sum _{j=u+1} ^n  r s_j, 
$$
equivalently, 
$$
(1-r) (\sum _{j=u+1}^n s_j ) = 1.
$$
This is impossible, because  $r$  and the all $s_j$'s are positive integers.
And we conclude that there is no such exact sequence as in  (14).
\qed
\end{pf}

%%%%%%%%%%%%%%%%%%%%%% References %%%%%%%%%%%%%%%%%%%%%%

%\end{article}


\begin{thebibliography}{}

\bibitem{AB}
{\sc M.Auslander} and {\sc M.Bridger},
{\em \lq\lq Stable module theory\rq\rq},
Mem. Amer. Math. Soc. 94 (1969).

\bibitem{AB2}
{\sc M.Auslander} and {\sc R.-O.Buchweitz},
{\em The homological theory of Cohen-Macaulay approximations},
 Mem. Soc. Math. de France {\bf 38} (1989), 5--37

\bibitem{AR}
{\sc M.Auslander} and {\sc I.Reiten}, 
{\em Applications of contravariantly finite subcategories} 
{\rm Adv. Math.} {\bf  86} (1991), 111--152.

\bibitem{AGP}
{\sc L.Avramov},  {\sc V.Gasharov},  and  {\sc I.Peeva},
{\em Complete Intersection dimension},
{\rm Inst. Hautes \'{E}tudes Sci. Publ. Math.} {\bf 86} (1997), 67-114 (1998).


\bibitem{Bass}
{\sc H.Bass},
{\em On the ubiquity of Gorenstein rings}
{\rm Math. Zeitschr.} {\bf 82} (1963), 8--28.


\bibitem{BH}
{\sc W.Bruns} and  {\sc J.Herzog},
{\em \lq\lq Cohen-Macaulay rings\rq\rq, revised version},
Cambridge University Press, 1998.


\bibitem{M}
{\sc H.Matsumura},
{\em \lq\lq Commutative ring theory\rq\rq},
Cambridge University Press, 1986.

\bibitem{V}
{\sc O.Veliche},
{\em Construction of modules with finite homological dimension}, 
{\rm J. Algebra}  {\bf 250} (2002), 427--449.

\bibitem{W}
{\sc T.Wakamatsu},
{\em Stable equivalence for self-injective algebras and a generalization of tilting modules}, 
{\rm J. Algebra} {\bf 134} (1990), 298--325.


\bibitem{X}
{\sc J.Xu},
{\em \lq\lq Flat covers of modules\rq\rq},
{\rm Lecture Notes in Math.}{\bf 1634}, 
Springer-Verlag, 1996.   

\bibitem{Y}
{\sc Y.Yoshino},
{\em \lq\lq Cohen-Macaulay modules over Cohen-Macaulay rings\rq\rq},
Lecture Note Sries {\bf 146}, Cambridge University Press, 1990. 

\end{thebibliography}
\end{document}